\documentclass[12pt]{article}
\usepackage{fullpage}
\usepackage{amsmath,amssymb,amsthm}
\usepackage{hyperref}
\hypersetup{
  colorlinks=true,
  linkcolor=blue,
  citecolor=blue,
  urlcolor=blue,
  pdftitle={Fishburn Numbers Count Homogeneous Games}
}
\newtheorem{theorem}{Theorem}
\newtheorem{lemma}[theorem]{Lemma}
\newtheorem{proposition}[theorem]{Proposition}
\newtheorem{corollary}[theorem]{Corollary}

\theoremstyle{remark}
\newtheorem{remark}[theorem]{Remark}

\newcommand{\seqnum}[1]{\href{https://oeis.org/#1}{\textup{\underline{#1}}}}

\title{Fishburn Numbers Count Homogeneous Games}
\author{M. Joris}
  
\date{}
\begin{document}

\maketitle

\begin{abstract}
We prove a new combinatorial interpretation for Fishburn numbers, in that they count the number of homogeneous games; specifically homogeneous games without dummy players, which we show extends to homogeneous games via the partial sums of Fishburn numbers. In doing so, we settle two conjectures commented in the OEIS entry \seqnum{A189359}, and finish off by providing some immediate consequences on homogeneous game counts; including an asymptotic formula and congruences.

\end{abstract}

\section{Introduction}
Homogeneous games are a class of weighted majority games, first introduced and studied by von Neumann and Morgenstern~\cite{vonNeumannMorgenstern}, and Shapley~\cite{Shapley1962} later studied them more generally. Their enumeration
has been studied by Sudh\"olter~\cite{Sudholter} and
Krohn and Sudh\"olter~\cite{KrohnSudholter}.
In this paper, we show that the numbers of homogeneous
games, up to the relabeling of the players, are counted by the
partial sums of the Fishburn numbers. We let $h_n$ denote the number of homogeneous games with $n$ players up to relabeling, with $h_0=0$, and let $\xi(n)$ denote the $n$-th Fishburn number. We write $F(1-q)$ for
the ordinary generating function of $(\xi)_{n\geq0}$.
These sequences are recorded in the OEIS as
\seqnum{A189359} and \seqnum{A022493}, respectively. In a comment dated December 13, 2021, on \seqnum{A189359},
Peter Bala~\cite{OEIS} made the astute observation that
$h_n-h_{n-1}=\xi(n)$ for $1\leq n\leq9$, and asked whether
this equality holds for all $n$. We answer this definitively.
We will start by stating our two main theorems. Note that we will introduce all of our notation precisely in Section~\ref{preliminaries}
\begin{theorem}\label{g_theorem}
Let the number of homogeneous games with no dummies be denoted $g_n$, then for every integer $n\geq1$,
\begin{equation}\label{g_fishburn_formula}
g_n=\xi(n).
\end{equation}
\end{theorem}
\begin{theorem}\label{main_theorem}
For every integer $n\geq1$,
\begin{equation}\label{h_difference_eq}
h_n-h_{n-1}=\xi(n).
\end{equation}
Equivalently,
\begin{equation}\label{h_sum_eq}
h_n=\sum_{j=1}^{n}\xi(j),
\end{equation}
and
\begin{equation}\label{h_generating_function}
\sum_{n\geq0}h_nq^n=\frac{F(1-q)-1}{1-q}.
\end{equation}
\end{theorem}
To prove Theorem~\ref{g_theorem}, we start from Sudh\"olter's recursive enumeration of incidence vectors associated with homogeneous games (see Section \ref{sudholter_enumeration}), and translate these recurrences into a functional equation
for a bivariate generating function that we solve iteratively. By comparing the resulting expression with
Zagier's generating function for the Fishburn numbers~\cite[Theorem~1]{Zagier}, we then prove Theorem \ref{g_theorem}. The formula for $h_n$ in Theorem \ref{main_theorem} then follows by deleting dummy players from the game, which we will make rigorous in Lemma \ref{deletion_extension_lemma}.

Regarding the paper's structure, we present the prerequisites in Section~\ref{preliminaries},
and review Sudh\"olter's enumeration and the recurrences
needed to start the proof. Section~\ref{generating_functions}
introduces the associated generating functions that we will need in Section~\ref{functional_eq_section} to derive and solve
the functional equation that proves Theorem \ref{g_theorem} and \ref{main_theorem}. We finish off by providing further comments and results in Section~\ref{consequences}.

\section{Preliminaries}\label{preliminaries}
We begin by the basic introduction of simple games to make sense of homogeneous games for readers unfamiliar.
\subsection{Homogeneous games}
A \emph{simple game} on the set of players $N=\{1,\ldots,n\}$ is a pair $(N,v)$,
where 
\begin{displaymath}
v:2^N\to\{0,1\}
\end{displaymath}
$v$ is called the \emph{characteristic function} and satisfies:
\begin{displaymath}
v(\varnothing)=0
\end{displaymath}
\begin{displaymath}
v(N)=1
\end{displaymath}
\begin{displaymath}
v(S)\leq v(T)
\end{displaymath}
whenever $S\subseteq T$. A \emph{coalition} is a subset of players $S\subseteq N$. It is winning if $v(S)=1$ and losing otherwise. A winning coalition is called \emph{minimal} if all its proper subsets are losing. 
A simple game is a \emph{weighted majority game} if there exist
non-negative integer weights $w_1,\ldots,w_n$ and an integer quota $\mathcal{Q}$,
with
\begin{displaymath}
0<\mathcal{Q}\leq\sum_{i=1}^{n}w_i
\end{displaymath}
such that a coalition $S$ is winning exactly when
$\sum_{i\in S}w_i\geq \mathcal{Q}$. The tuple $(\mathcal{Q};w_1,\ldots,w_n)$ is called a \emph{representation} of the game.
The game is \emph{homogeneous} if the weights and quota can be chosen
so that every minimal winning coalition has total weight exactly $\mathcal{Q}$. We use the general formulation of Sudh\"olter~\cite{Sudholter},
which allows games in which a coalition and its complement both be winning or both be losing.

\subsection{Counting homogeneous games, and our notation}
We let $h_n$ denote the number of homogeneous games on $n$ players up to relabeling of the players, and set $h_0=0$.
Two games are counted as \emph{equivalent} if there is a bijection between their player sets that preserves winning and losing coalitions; different weighted representations of the same game are not counted separately (or else there would be infinitely many and this counting problem would be trivial). For further reading on the topic, see Peleg and
Sudh\"olter~\cite[Sections~2.2.3 and~2.3]{PelegSudholter}.
Krohn and Sudh\"olter~\cite[Table~1]{KrohnSudholter} enumerate $h_n$ up to nine players, and the initial values are

\begin{displaymath}
0,1,3,8,23,76,293,1307,6642,37882
\end{displaymath}
as recorded in the OEIS~\cite{OEIS} entry \seqnum{A189359}. Furthermore, a player $i$ is called a \emph{dummy} if $S$ and $S\cup\{i\}$
have the same winning or losing status for every
$S\subseteq N\setminus\{i\}$.

We let $g_n$ denote the number of homogeneous games on $n$ players
without dummies, counted up to relabeling.

We let $b_n$ count those games that have no dummies and no one-person winning coalitions.
In particular, $b_1=0$, since the player in a one-player game forms a winning coalition.

\subsection{Fishburn numbers}
Let $(a;q)_n := (1-a)(1-aq)\cdots(1-aq^{n-1})$ denote the $q$-Pochhammer symbol, with $(q)_n := (q;q)_n$, then the following series
\begin{displaymath}
    F(q) = \sum_{n\geq 0} (q)_n,
\end{displaymath}
is the Kontsevich-Zagier series commonly monikered the ``strange'' series \cite{Zagier}. Fishburn numbers \seqnum{A022493}, denoted $\xi(n)$, are given as the coefficients of $F(1-q)$,
\begin{displaymath}
F(1-q) = \sum_{n\geq0}\xi(n)q^n = 1+q+2q^2+5q^3+15q^4+53q^5+217q^6+\ldots
\end{displaymath}
Zagier~\cite[Theorem 1]{Zagier} gives the generating function as
\begin{equation}\label{fishburn_generating_funct}
F(1-q)=\sum_{n\geq0}\xi(n)q^n = 
\sum_{n\geq 0}(1-q)_n=\sum_{m\geq0}\prod_{i=1}^{m}\bigl(1-(1-q)^i\bigr).
\end{equation}
Bousquet-M\'elou et al.\ \cite{BMCDK} relate Fishburn numbers to several classes of combinatorial objects; including ascent sequences and unlabeled $(2+2)$-free posets. These are but a few of many objects that Fishburn numbers count.
Furthermore, there is a growing body of work relating $F(q)$ and generalizations thereof to modular forms. For some further reading see \cite{GuerzhoyKentRolen2014}, \cite{Zagier2010}, \cite{AhlgrenKim2015}.

\subsection{Sudh\"olter's enumeration}\label{sudholter_enumeration}

Sudh\"olter~\cite[Lemma~1.1]{Sudholter} reduces the
enumeration of homogeneous games by deleting dummy
players and then players who win alone. We give the following lemma to glue together certain arguments that we will require.
\begin{lemma}\label{deletion_extension_lemma}
For a homogeneous game with $n$ players:
\begin{enumerate}
\item Deleting all dummy players leaves a well-defined
homogeneous game without dummies.

\item Adjoining dummy players to restore $n$ players is inverse
to this operation, up to relabeling.

\item With no dummies, and not every nonempty coalition
wins, deleting all players who win alone leaves a well-defined
homogeneous game with neither dummies nor singleton winning
coalitions.

\item Adjoining players who win alone to restore $n$ players
is inverse to the operation in (3), up to relabeling.
\end{enumerate}
\end{lemma}
\begin{proof}
These are due to \cite[Lemma~1.1]{Sudholter}. For (1) and (2), let $D$ be the set of dummy players, put $R=N\setminus D$, and define $v_R(T):=v(T)$ for $T\subseteq R$. Since deleting dummy players does not
change a coalition's value we get $v(S)=v_R(S\cap R)$ for $S\subseteq N$. Thus $v_R(R)=v(N)=1$, and any dummy in $v_R$ would also
be a dummy in $v$. Conversely, the displayed formula
defines the unique extension of $v_R$ by a given
set $D$ of dummy players. If $v_R$ has no dummies,
the dummy set of this extension is exactly $D$.
Hence, the key idea is that deletion and extension are inverse.

For (3) and (4), suppose that $v$ has no dummies.
Let $W=\{i\in N:v(\{i\})=1\}$ with $R=N\setminus W$. If $R=\varnothing$, every player wins alone. For any nonempty
coalition $S$, choose $i\in S$. Since $\{i\}\subseteq S$, the
defining condition for simple games gives
$1=v(\{i\})\leq v(S)$, so $S$ wins. Thus we get the unique game
in which every nonempty coalition wins.

Suppose now that $R\ne\varnothing$. If $v(R)=0$, then any
coalition $S$ containing a player $j\in W$ wins, since
$\{j\}\subseteq S$ and $v(\{j\})=1$. Any coalition disjoint from $W$ is a subset of $R$, so it loses since $v(S)\leq v(R)=0$. Thus the winning coalitions have at least one player in $W$, making every player in $R$ a dummy. Therefore $v(R)=1$, and restriction to $R$ defines a simple game $v_R$.
It has no singleton winning coalitions by definition of $W$.

For each $i\in R$, choose $T\subseteq N\setminus\{i\}$ such that $v(T)=0$ and $v(T\cup\{i\})=1$. Such a coalition exists because $i$ is not a dummy. Since $T$ is losing, it contains no player from $W$: otherwise, for some $j\in T\cap W$, we would have $\{j\}\subseteq T$ and hence $1=v(\{j\})\leq v(T)$. Consequently, the exceptional game above is the only case we need to consider separately in which restriction to $R$ fails to give the required simple game. The original game is recovered by
\[
v(S)=
\begin{cases}
1, & S\cap W\ne\varnothing,\\
v_R(S\cap R), & S\cap W=\varnothing.
\end{cases}
\]
Conversely, this formula defines the unique extension by a prescribed set $W$ of players who win alone.
Its singleton winners are exactly $W$, and the players
of $R$ remain non-dummies. Hence these operations are
also inverse. Both restrictions preserve homogeneity, since every minimal winning coalition of a restricted game is also
minimal winning in the original game. For the extensions,
take any homogeneous representation of the restricted
game and assign each added player a weight zero in (2),
or weight equal to the quota in (4).

Finally, the deleted sets are determined by the
characteristic function. Any isomorphism between
restricted games extends by any bijection between
their respective sets of added players.
\end{proof}
In our notation, these reductions give, for $n\geq1$,
\begin{align}
h_n&=\sum_{j=1}^{n}g_j,
\label{dummy_deletion}\\
g_n&=1+\sum_{j=2}^{n}b_j,
\label{singleton_deletion}
\end{align}
where the additional `$1$' term in \eqref{singleton_deletion} counts
the unique game in which every nonempty coalition wins.
For the first identity, deleting all dummy players leaves
a game without dummies on some $j$ players. Conversely,
each such game extends uniquely, up to relabeling, to a
game on $n$ players by then adjoining $n-j$ dummy players.

We want to make explicit that Lemma \ref{deletion_extension_lemma} belongs to Sudh\"olter, and that we merely state his ideas with a refined proof for the sake of this paper's self-containedness.

To enumerate these remaining games, Sudh\"olter uses vectors he calls \emph{incidence vectors} and
obtains a recursive formula for them. We use
the resulting recurrences, which we state below in
terms of the coefficients $a_n^k$.

\begin{proposition}[Sudh\"olter]\label{sudholter_enumeration_prop}
Let $(a_n^k)_{n\geq2,\,k\geq1}$ be defined by
the initial values
\begin{equation}\label{initial_values}
a_2^1=1,
\qquad
a_3^1=0,
\end{equation}
the convention $a_n^k=0$ when $n<2k$, and the following
recurrences. For $n\geq3$,
\begin{equation}\label{one_recurrence_eq}
a_{n+1}^1
=\sum_{j=1}^{\lfloor n/2\rfloor}
 \sum_{p=0}^{n-2j}
 \binom{j+p-1}{p}a_{n-p}^{j}.
\end{equation}
For $k\geq2$ and $n\geq2k-1$,
\begin{equation}\label{multiple_recurrence_eq}
a_{n+1}^{k}
=\sum_{j=k-1}^{\lfloor n/2\rfloor}
 \sum_{p=0}^{n-2j}
 \binom{j}{k-1}\binom{j+p-1}{p}a_{n-p}^{j}
 -a_n^{k-1}.
\end{equation}
Then, for every $n\geq2$,
\begin{equation}\label{reduced_game_count_eq}
b_n=a_{n+2}^1.
\end{equation}
\end{proposition}
The recurrences are Theorem~3.12(i) and~(ii)
of~\cite{Sudholter}, with the sums reordered in the
second recurrence. The value $b_0 = a_2^1=1$ is auxiliary, as in Sudh\"olter's original enumeration. The relation \eqref{reduced_game_count_eq} combines
Sudh\"olter's Lemma~1.4 and Theorem~2.10 with
Theorem~3.8(i). The first two results identify the games
counted by $b_n$, and the last result states that this number is $a_{n+2}^1$. These identities provide the counting relations used in the generating function argument below.

\section{Generating functions}\label{generating_functions}

Define
\begin{displaymath}
G(q)=\sum_{n\geq1}g_nq^n,
\qquad
B(q)=\sum_{n\geq2}b_nq^n.
\end{displaymath}
The relation \eqref{singleton_deletion} between their coefficients gives
\begin{align}
G(q) =\sum_{n\geq1}\left(1+\sum_{j=2}^{n}b_j\right)q^n =\frac{q}{1-q}
  +\sum_{j\geq2}b_j\sum_{n\geq j}q^n =\frac{q+B(q)}{1-q}.
\label{g_b_eq}
\end{align}
We next introduce generating functions for the coefficients $a_n^k$,
first keeping $k$ fixed. For each $k\geq1$, define
\begin{align}
A_k(q)=\sum_{n\geq2k}a_n^kq^n,\\
C_k(q)=\frac{A_k(q)}{(1-q)^k}.
\end{align}
The introduction of $C_k(q)$ makes the calculation in our later Lemmas cleaner. The factor $(1-q)^{-k}$ is chosen to account for the
binomial coefficients in the recurrences as, by the binomial theorem,
\begin{displaymath}
(1-q)^{-k} =\sum_{p\geq0}\binom{k+p-1}{p}q^p.
\end{displaymath}
Thus multiplying $A_k(q)$ by $(1-q)^{-k}$ incorporates
the summation over $p$ in \eqref{one_recurrence_eq}
and \eqref{multiple_recurrence_eq}.
To keep track of the dependence on $k$, we introduce a second variable $z$
and set
\begin{equation}\label{bivariate_series_eq}
f(z)=\sum_{k\geq1}C_k(q)z^k,
\end{equation}
where we suppress the dependence on $q$ in this notation.

\begin{remark}\label{formal_substitutions}
Since $C_k(q)\in q^{2k}\mathbb{Q}[[q]]$, the coefficient of
$q^n$ in $f(z)$ is a polynomial in $z$ of degree at most
$\lfloor n/2\rfloor$. Thus $f(z)\in\mathbb{Q}[z][[q]]$.
For every $\phi(q)\in\mathbb{Q}[[q]]$, the substitution
\begin{displaymath}
f(\phi(q))=\sum_{k\geq1}C_k(q)\phi(q)^k,
\end{displaymath}
is well defined in $\mathbb{Q}[[q]]$, since only terms with
$2k\leq n$ can contribute to the coefficient of $q^n$.
\end{remark}

\begin{lemma}\label{normalization_lemma}
We have
\begin{equation}\label{normalization_eq}
G(q)
=\sum_{k\geq1}\frac{A_k(q)}{q(1-q)^{k+1}}
=\frac{f(1)}{q(1-q)}.
\end{equation}
\end{lemma}

\begin{proof}
The substitution $z=1$ is valid by
Remark~\ref{formal_substitutions}. Expanding $(1-q)^{-k}$ by the Binomial theorem and collecting terms of equal degree gives
\begin{align*}
f(1) =\sum_{k\geq1}\frac{A_k(q)}{(1-q)^k}=\sum_{k\geq1}
  \left(\sum_{m\geq2k}a_m^kq^m\right)
  \left(\sum_{p\geq0}\binom{k+p-1}{p}q^p\right)
\end{align*}
Only finitely many terms contribute to each coefficient. Next, put $n=m+p$, then
\begin{displaymath}
f(1) =\sum_{k\geq1}\left(\sum_{m\geq2k}a_m^kq^m\right)\left(\sum_{p\geq0}\binom{k+p-1}{p}q^p\right)=\sum_{n\geq2}\left(\sum_{k=1}^{\lfloor n/2\rfloor}\sum_{p=0}^{n-2k}\binom{k+p-1}{p}a_{n-p}^k\right)q^n
\end{displaymath}
For $n\geq3$, the recurrence \eqref{one_recurrence_eq}
and the relation \eqref{reduced_game_count_eq} give
\begin{displaymath}
\sum_{k=1}^{\lfloor n/2\rfloor}\sum_{p=0}^{n-2k}\binom{k+p-1}{p}a_{n-p}^k=a_{n+1}^1=b_{n-1}.
\end{displaymath}
For $n=2$, the sum consists only of the auxiliary
term $a_2^1=1$. Hence,
\begin{equation}\label{f_one_eq}
f(1)=q^2+\sum_{n\geq3}b_{n-1}q^n=q^2+qB(q).
\end{equation}
Combining this identity with \eqref{g_b_eq} yields
\begin{displaymath}
G(q) =\frac{q+B(q)}{1-q} =\frac{f(1)}{q(1-q)} =\sum_{k\geq1}\frac{A_k(q)}{q(1-q)^{k+1}},
\end{displaymath}
as required.
\end{proof}

\section{Proof of the main theorem}\label{functional_eq_section}
\begin{lemma}\label{functional_eq_lemma}
The series $f$ satisfies
\begin{equation}\label{functional_eq}
(1+qz)f((1-q)z)=qzf(1+z)+q^2(1-q)z.
\end{equation}
\end{lemma}

\begin{proof}
Using $a_2^1=1$, $a_3^1=0$, and
\eqref{reduced_game_count_eq}, we have
\begin{align*}
A_1(q)
&=q^2+\sum_{n\geq2}b_nq^{n+2} =q^2+q^2B(q) =q^2+q\bigl(f(1)-q^2\bigr),
\end{align*}
where the last equality follows from \eqref{f_one_eq}.
Hence
\begin{equation}\label{first_series}
(1-q)C_1(q)=q^2(1-q)+qf(1).
\end{equation}
Fix $k\geq2$, then recurrence~\eqref{multiple_recurrence_eq} reads
\begin{displaymath}
a_{n+1}^{k} = \sum_{j=k-1}^{\lfloor n/2\rfloor}\sum_{p=0}^{n-2j}\binom{j}{k-1}\binom{j+p-1}{p}a_{n-p}^{j}-a_n^{k-1}.
\end{displaymath}
for $n\geq2k-1$. At $n=2k-2$, the double sum in
\eqref{multiple_recurrence_eq} consists only of the term
$a_{2k-2}^{k-1}$, which cancels the subtracted term.
The left-hand side is $a_{2k-1}^{k}=0$ by convention.
Thus the recurrence also holds at $n=2k-2$. Multiplying by $q^{n+1}$, summing over $n\geq2k-2$,
and changing the order of summation with $m=n-p$ gives
\begin{align*}
A_k(q)&=q\sum_{j\geq k-1}\binom{j}{k-1}A_j(q)\sum_{p\geq0}\binom{j+p-1}{p}q^p-qA_{k-1}(q)\\
&=q\sum_{j\geq k-1}\binom{j}{k-1}C_j(q)-qA_{k-1}(q).
\end{align*}
Substituting in $A_j(q)=(1-q)^jC_j(q)$ gives
\begin{equation}\label{remaining_series_eq}
(1-q)^kC_k(q)=q\sum_{j\geq k-1}\binom{j}{k-1}C_j(q)-q(1-q)^{k-1}C_{k-1}(q)
\end{equation}
Multiplying \eqref{remaining_series_eq} by $z^k$ and summing
over $k\geq2$ gives
\begin{displaymath}
\sum_{k\geq2}(1-q)^kC_k(q)z^k =q\sum_{k\geq2}z^k\sum_{j\geq k-1}\binom{j}{k-1}C_j(q)-q\sum_{k\geq2}(1-q)^{k-1}C_{k-1}(q)z^k.
\end{displaymath}
In the first term on the right, the conditions
$k\geq2$ and $j\geq k-1$ are equivalent to
$j\geq1$ and $2\leq k\leq j+1$. For each $j\geq1$, the binomial theorem, with $r=k-1$, gives
\begin{displaymath}
\sum_{k=2}^{j+1}\binom{j}{k-1}z^k=z\sum_{r=1}^{j}\binom{j}{r}z^r=z\bigl((1+z)^j-1\bigr).
\end{displaymath}
We can therefore reverse the order of summation, put $r=k-1$,
and apply the binomial theorem as above to get
\begin{align*}
q\sum_{k\geq2}z^k\sum_{j\geq k-1}\binom{j}{k-1}C_j(q)
&=q\sum_{j\geq1}C_j(q)\sum_{k=2}^{j+1}\binom{j}{k-1}z^k\\
&=qz\sum_{j\geq1}C_j(q)\sum_{r=1}^{j}\binom{j}{r}z^r\\
&=qz\sum_{j\geq1}C_j(q)\bigl((1+z)^j-1\bigr)\\
&=qz\bigl(f(1+z)-f(1)\bigr).
\end{align*}
Note again that these rearrangements are valid coefficientwise because of Remark \ref{formal_substitutions}. The summed recurrence is therefore
\begin{displaymath}
f((1-q)z)-(1-q)C_1(q)z=qz\bigl(f(1+z)-f(1)\bigr)-qzf((1-q)z).
\end{displaymath}
Substituting in \eqref{first_series} gives
\begin{displaymath}
f((1-q)z)-(q^2(1-q)+qf(1))z=qz\bigl(f(1+z)-f(1)\bigr)-qzf((1-q)z).
\end{displaymath}
and trivial rearranging yields \eqref{functional_eq} as desired.
\end{proof}

\begin{lemma}\label{solution_lemma}
We have
\begin{equation}\label{series_solution_eq}
f(1)=q(1-q)\sum_{m\geq1}\prod_{i=1}^{m}\bigl(1-(1-q)^i\bigr) = q(1-q)(F(1-q)-1).
\end{equation}
\end{lemma}

\begin{proof}
A very natural way to prove this is by the substitution of a recursion $u\in\mathbb{Q}[[q]]$ that brings out the $1-(1-q)^i$ factors whilst giving us $f(1)$. Indeed, this is what we will do. Substituting $z=u/(1-q)$
into \eqref{functional_eq} and multiplying by $1-q$ gives us
\begin{displaymath}
(1-q+qu)f(u)=qu f\left(1+\frac{u}{1-q}\right)+q^2(1-q)u.
\end{displaymath}
Since $(1-q+qu)$ has constant term $1$, it is invertible
in $\mathbb{Q}[[q]]$. Hence
\begin{equation}\label{substitution_eq}
f(u)=\frac{qu}{1-q+qu}f\left(1+\frac{u}{1-q}\right)+\frac{q^2(1-q)u}{1-q+qu}.
\end{equation}
Now we let $u_0=1$, and define for $r\geq0$
\begin{equation}\label{u_recurrence}
u_{r+1}=1+\frac{u_r}{1-q}.
\end{equation}
Induction on $r$ gives
\begin{equation}\label{iteration_eq}
u_r=\sum_{j=0}^{r}(1-q)^{-j} =(1-q)^{-r}\sum_{j=0}^{r}(1-q)^j =\frac{1-(1-q)^{r+1}}{q(1-q)^r}.
\end{equation}
In particular, each $u_r$ belongs to $\mathbb{Q}[[q]]$, so $f(u_r)$ is well defined, again by Remark~\ref{formal_substitutions}. Moreover, this now provides the $1-(1-q)^i$ factors as mentioned, so we can go on simplifying \eqref{substitution_eq} by writing 
\begin{displaymath}
1-q+qu_r =1-q+\frac{1-(1-q)^{r+1}}{(1-q)^r} =(1-q)^{-r},
\end{displaymath}
and
\begin{displaymath}
\frac{qu_r}{1-q+qu_r} =q(1-q)^ru_r =1-(1-q)^{r+1},
\end{displaymath}
and then substituting these expressions into \eqref{substitution_eq} and using \eqref{u_recurrence}, we obtain
\begin{equation}\label{iteration_step}
f(u_r)=\bigl(1-(1-q)^{r+1}\bigr)f(u_{r+1}) +q(1-q)\bigl(1-(1-q)^{r+1}\bigr).
\end{equation}
We briefly remind the reader of the shorthand notation for the $q$-Pochhammer symbol $(q)_n := (q;q)_n$ with $(1-q)_0=1$, which allows us to simplify the proof visually. Then for $m\geq1$ and $r$ as the summation variable, then since $(1-q)_{r+1}=(1-q)_r\bigl(1-(1-q)^{r+1}\bigr)$,
multiplying \eqref{iteration_step} by $(1-q)_r$ and rearranging gives
\begin{displaymath}
(1-q)_rf(u_r)-(1-q)_{r+1}f(u_{r+1})=q(1-q)(1-q)_{r+1}.
\end{displaymath}
For $N\geq1$ and summing over $r=0,\ldots,N-1$, the left-hand side telescopes. Since $(1-q)_0=1$ and $u_0=1$, this yields
\begin{equation}\label{finite_iteration_eq}
f(1)=q(1-q)\sum_{m=1}^{N}(1-q)_m+(1-q)_Nf(u_N).
\end{equation}
Each factor $1-(1-q)^i$ is divisible by $q$, so $(1-q)_m\in q^m\mathbb{Q}[[q]]$. Thus $\sum_{m\geq1}(1-q)_m$ is well defined as a formal power
series. Also, with Remark \ref{formal_substitutions} we have
\begin{displaymath}
f(u_N)=\sum_{k\geq1}C_k(q)u_N^k\in q^2\mathbb{Q}[[q]].
\end{displaymath}
Consequently, $(1-q)_Nf(u_N)\in q^{N+2}\mathbb{Q}[[q]]$, so the remainder tends to zero coefficientwise as $N\to\infty$. Taking the coefficientwise limit in \eqref{finite_iteration_eq} proves \eqref{series_solution_eq}.
\end{proof} 
\begin{proof}[Proof of Theorem~\ref{g_theorem} and Theorem~\ref{main_theorem}]
Lemma~\ref{normalization_lemma} and Lemma~\ref{solution_lemma},
together with \eqref{fishburn_generating_funct}, give
\begin{align*}
G(q) =\frac{f(1)}{q(1-q)} =\sum_{m\geq1}\prod_{i=1}^{m}\bigl(1-(1-q)^i\bigr) =F(1-q)-1.
\end{align*}
Comparing coefficients yields $g_n=\xi(n)$, proving Theorem \ref{g_theorem}. For every $n\geq1$.
By \eqref{dummy_deletion}, we therefore have
\begin{displaymath}
h_n=\sum_{j=1}^{n}g_j=\sum_{j=1}^{n}\xi(j),
\end{displaymath}
which proves \eqref{h_sum_eq}. Taking successive differences
and using $h_0=0$ gives $h_n-h_{n-1}=\xi(n)$ for $n\geq1$ as asserted in \eqref{h_difference_eq}. Finally, changing the order of summation gives
\begin{align*}
\sum_{n\geq0}h_nq^n =\sum_{n\geq1}\sum_{j=1}^{n}\xi(j)q^n\ =\sum_{j\geq1}\xi(j)\sum_{n\geq j}q^n =\frac{1}{1-q}\sum_{j\geq1}\xi(j)q^j\ =\frac{F(1-q)-1}{1-q},
\end{align*}
which proves \eqref{h_generating_function}, all-together proving Theorem \ref{main_theorem}.
\end{proof}

\section{Further results}\label{consequences}

Any formulas involving $\xi(n)$ now supply formulas on $g_n$, and consequently $h_n$. Whilst the purpose of this paper is to establish the Fishburn-Homogeneous games connection, we will provide some immediate consequences on $g_n$ and $h_n$ of the well-known properties of Fishburn numbers from the literature. This includes an asymptotic formula for $h_n$ (Corollary \ref{asymptotic_growth_corollary}), congruences on $h_n$ (Corollary \ref{congruences_corollary}), and a simple probabilistic result (Corollary \ref{probability_corollary}). We also provide a solution to Gladkovskii's conjecture commented in the OEIS entry \seqnum{A189359} (Corollary \ref{gladkovskii_conjecture}). We then finish off by remarking that bijections of homogeneous games with other objects also counted by Fishburn numbers must exist.

\begin{corollary}
\label{gladkovskii_conjecture}
In a comment dated May 13 2013 in \seqnum{A189359}, Sergei Gladkovskii conjectured that for $k\geq0$, and letting $Q_k(q)\in1+q\mathbb{Q}[[q]]$ be given by the continued fraction
\begin{displaymath}
Q_k(q)=1+\frac{1-(1-q)^{2k+2}}{1-\dfrac{1-(1-q)^{2k+3}}{1-(1-q)^{2k+3}+1/Q_{k+1}(q)}},
\end{displaymath}
then
\begin{displaymath}
\sum_{n\geq0}h_nq^n=\frac{qQ_0(q)}{1-q}.
\end{displaymath}
\end{corollary}

\begin{proof}
We can simplify,
\begin{displaymath}
Q_k(q) =1+\bigl(1-(1-q)^{2k+2}\bigr)+\bigl(1-(1-q)^{2k+2}\bigr)\bigl(1-(1-q)^{2k+3}\bigr)Q_{k+1}(q).
\end{displaymath}
Notice,
\begin{displaymath}
Q_k(q) =1+q\sum_{i=0}^{2k+1}(1-q)^i+q^2\left(\sum_{i=0}^{2k+1}(1-q)^i\right)\left(\sum_{j=0}^{2k+2}(1-q)^j\right)Q_{k+1}(q),
\end{displaymath}
using geometric series identity for $1-(1-q)^m$. The following can be deduced from rather tedious successive substitution and comparing against its previous substitution,
\begin{displaymath}
Q_0(q) =\frac{1}{q}\sum_{m=0}^{2N-1}(1-q)_{m+1} +\frac{(1-q)_{2N+1}}{q}Q_N(q).
\end{displaymath}
we can now use the same argument as we used at the end of the proof for Lemma \ref{solution_lemma} for the vanishing of the remainder terms, to find that
\begin{displaymath}
Q_0(q)=\sum_{m\geq0}\prod_{j=2}^{m+1}\bigl(1-(1-q)^j\bigr).
\end{displaymath}
The rest follows by Theorem~\ref{main_theorem}, and we find that
\begin{displaymath}
    qQ_0(q)=F(1-q)-1.
\end{displaymath}
proving the conjecture.
\end{proof}

\begin{corollary}\label{asymptotic_growth_corollary}
As $n\to\infty$, we have the following asymptotic formula
\begin{equation}\label{asymptotic_growth_eq}
h_n \sim \frac{12\sqrt{3}}{\pi^{5/2}}e^{\pi^2/12}\left(\frac{6}{\pi^2}\right)^n n!\sqrt{n}.
\end{equation}
Moreover,
\begin{equation}\label{ratio_limit_eq}
\lim_{n\to\infty}\frac{h_{n+1}}{nh_n} =\frac{6}{\pi^2}.
\end{equation}
\end{corollary}

\begin{proof}
Zagier~\cite[Theorem~4]{Zagier} proves that
\begin{displaymath}
\xi(n)\sim
\frac{12\sqrt{3}}{\pi^{5/2}}e^{\pi^2/12}\left(\frac{6}{\pi^2}\right)^n n!\sqrt{n}.
\end{displaymath}
To obtain \eqref{asymptotic_growth_eq}, it remains to show that $h_n\sim \xi(n)$. The asymptotic formula gives
\begin{displaymath}
\frac{\xi(n-1)}{\xi(n)} \sim \frac{\pi^2}{6n}\longrightarrow0.
\end{displaymath}
Thus $(\xi(n))_n$ is eventually strictly increasing and tends
to infinity. By the Stolz--Ces\`aro theorem and
\eqref{h_difference_eq},
\begin{displaymath}
\lim_{n\to\infty}\frac{h_n}{\xi(n)} = \lim_{n\to\infty}\frac{h_n-h_{n-1}}{\xi(n)-\xi(n-1)} = \lim_{n\to\infty}\frac{\xi(n)}{\xi(n)-\xi(n-1)} = 1,
\end{displaymath}
which proves \eqref{asymptotic_growth_eq}.
Finally, using $h_n\sim \xi(n)$ and the same asymptotic
formula, we obtain
\begin{displaymath}
\frac{h_{n+1}}{nh_n}\sim\frac{\xi(n+1)}{n\xi(n)}\longrightarrow\frac{6}{\pi^2} ,
\end{displaymath}
proving \eqref{ratio_limit_eq}.
\end{proof}
\begin{corollary}\label{congruences_corollary}
For every integer $n\geq0$,
\begin{equation}\label{congruence_modulo_five_eq}
h_{5n+2}\equiv h_{5n+3}\equiv h_{5n+4}\pmod{5}.
\end{equation}
More generally, for all integers $r\geq1$ and $n\geq1$,
\begin{equation}\label{congruence_modulo_powers_eq}
h_{5^r n-3}\equiv h_{5^r n-2}
\equiv h_{5^r n-1}\pmod{5^r}.
\end{equation}
\end{corollary}

\begin{proof}
Andrews and Sellers~\cite{AndrewsSellers} prove that
\begin{displaymath}
\xi(5n+3)\equiv \xi(5n+4)\equiv0\pmod{5}.
\end{displaymath}
for $n\geq0$. Applying \eqref{h_difference_eq} at $5n+3$ and $5n+4$
gives \eqref{congruence_modulo_five_eq}. Straub~\cite[Example~1.3]{Straub} extends these congruences to
\begin{displaymath}
\xi(5^r n-2)\equiv \xi(5^r n-1)\equiv0\pmod{5^r}.
\end{displaymath}
for $r\geq1,\ n\geq1$. Applying \eqref{h_difference_eq} at $5^r n-2$ and
$5^r n-1$ gives \eqref{congruence_modulo_powers_eq}.
\end{proof}
\begin{corollary}\label{probability_corollary}
Choose some homogeneous game uniformly from the set of all homogeneous games on $n$ players up to equivalence. Let $D$ denote the set of dummies in the chosen homogeneous game with $n$ players, then the probability that that game contains a dummy player is
\begin{equation}\label{probability_eq1}
P(D \neq \emptyset) = \frac{h_n-g_n}{h_n} =\frac{h_{n-1}}{h_n} \sim \frac{\pi^2}{6n}.
\end{equation}
Furthermore, choose uniformly some $i \in N$. Then,
\begin{equation}\label{probability_eq2}
P(i \in D) \sim \frac{\pi^2}{6n^2}.
\end{equation}
\end{corollary}

\begin{proof}
By \eqref{dummy_deletion}, the number of homogeneous
games on $n$ players with at least one dummy player is $h_n-g_n=\sum_{j=1}^{n-1}g_j=h_{n-1}$.
The asymptotic formula \eqref{asymptotic_growth_eq} gives
\begin{displaymath}
\frac{h_{n-1}}{h_n} \sim \frac{\pi^2}{6}\frac{(n-1)!\sqrt{n-1}}{n!\sqrt{n}} = \frac{\pi^2}{6n}\sqrt{1-\frac1n} \sim \frac{\pi^2}{6n}.
\end{displaymath}
which proves \eqref{probability_eq1}. For $P(i \in D)$ choose
a player $i$ uniformly. Deleting $d$ dummy players gives $P(\lvert D \rvert \geq d)=h_{n-d}/h_n$ for $1\leq d<n$. Hence, by the tail-sum formula,
\begin{displaymath}
P(i\in D)=\frac{\mathbb{E}[\lvert D \rvert]}{n}
=\frac{1}{nh_n}\sum_{j=1}^{n-1}h_j.
\end{displaymath}
Stolz--Ces\`aro gives
$\sum_{j=1}^{n-1}h_j\sim h_{n-1}$, and therefore
\[
P(i\in D)\sim\frac{h_{n-1}}{nh_n}
\sim\frac{\pi^2}{6n^2},
\]
proving \eqref{probability_eq2}.

\end{proof}

\begin{remark}
Since $g_n = \xi(n)$, any object that is counted by Fishburn numbers has some existing bijection with no-dummy homogeneous games, and extend more generally to homogeneous games by Lemma \ref{deletion_extension_lemma}.
\end{remark}

\section{Acknowledgments}\label{acknowledgments}
We would like to thank Sascha Kurz and Xavier Molinero for their helpful feedback on earlier drafts, helping improve the overall paper and preliminary game theory presented. We initially used GPT-5.4 for literature review and for understanding the necessary objects from Sudh\"olter~\cite{Sudholter} to construct the generating functions, and later a small correction to the limits in our original proof of Lemma \ref{functional_eq_lemma}.

2020 \emph{Mathematics Subject Classification}: Primary 05A15;
Secondary 05A16, 91A12.

\emph{Keywords}: Fishburn number, homogeneous game, weighted majority game, generating function.

(Concerned with sequences \seqnum{A022493} and \seqnum{A189359}.)
\end{document}